 \documentclass[a4paper,10pt]{article}
 
 \usepackage{xy,amssymb,amsmath,diagrams,times}
 
 \xyoption{all}
 
 \newcommand{\N}{\mathbb{N}}
 \newcommand{\C}{\mathbb{C}}
 \newcommand{\Z}{\mathbb{Z}}
 \newcommand{\wis}[1]{{\text{\em \usefont{OT1}{cmtt}{m}{n} #1}}}
 \newcommand{\Oscr}{\mathcal{O}}

\newcommand{\vtx}[1]{*+[o][F-]{\scriptscriptstyle #1}}

 \newtheorem{definition}{Definition}
 \newtheorem{proposition}{Proposition}
 
 \newtheorem{theorem}{Theorem}
 \newtheorem{lemma}{Lemma}
 \newtheorem{example}{Example}

  \newcommand{\proofhead}[1]{\par\pagebreak[1]\noindent{\em#1.\ }}
 \newcommand{\pf}{\proofhead{Proof}}
  \newcommand{\qed}{{\unskip\nolinebreak[1]\hspace{1.5em}\mbox{}\nolinebreak
    \hfill$\square$\parfillskip=0pt\finalhyphendemerits=0\par\pagebreak[1]}}

\newenvironment{proof}{\pf}{\qed \par \vskip 6mm}
 
 \title{Qurves and Quivers}
 \author{Lieven Le Bruyn \\ Department of Mathematics, University of Antwerp \\ 
 Middelheimlaan 1, B-2020 Antwerp (Belgium) \\ {\tt lieven.lebruyn@ua.ac.be}}
 \date{\empty}
 
\begin{document}
 \sloppy
 
 \mathversion{bold}
 
 \maketitle

 \begin{abstract}
 In this paper we associate to an {\em $\overline{\ell}$-qurve} $A$ (formerly known as a quasi-free algebra \cite{CuntzQuillen} or formally smooth algebra \cite{KontRos}) the {\em one-quiver} $Q_1(A)$ and dimension vector $\alpha_1(A)$. This pair contains enough information to reconstruct for all $n \in \N$ the $GL_n$-\'etale local structure of the representation scheme $rep_n~A$. In an appendix we indicate how one might extend this to qurves over non-algebraically closed fields.
 Further, we classify all finitely generated groups $G$ such that the group algebra $\ell G$ is an $\ell$-qurve. If $char(\ell) = 0$ these are exactly the virtually free groups. We determine the one-quiver setting in this case and indicate how it can be used to study the finite dimensional representations of virtually free groups. As this approach also applies to {\em fundamental algebras} of {\em graphs of separable $\ell$-algebras}, we state the results in this more general setting.
 \end{abstract}

 \section{Qurves}
 
 In this paper, $\ell$ is a commutative field with algebraic closure $\overline{\ell}$. Algebras will be associative $\ell$-algebras with unit and (usually) finitely generated over $\ell$.
 For an $\ell$-algebra $A$ let $A'$ be the $\ell$-vectorspace $A/\ell.1$ and define (following \cite[\S 1]{CuntzQuillen}) the graded algebra of {\em non-commutative differential forms}
 \[
 \Omega A = \oplus_{i=o}^{\infty}~\Omega^i A \qquad \text{with} \qquad \Omega^i~A = A \otimes A'^{\otimes i} \]
 with multiplication defined by the maps $\Omega^n A \otimes \Omega^{k-1} A \rTo \Omega^{n+k-1} A$ where
 \[
 (a_0,\hdots,a_n).(a_{n+1},\hdots,a_{n+k}) = \sum_{i=0}^n (-1)^{n-i} (a_0,\hdots,a_ia_{i+1},\hdots,a_{n+k}) \]
 As $\Omega^0 A = A$ this multiplication defines an $A$-bimodule structure on each $\Omega^n A$ and one proves \cite[Prop. 2.3]{CuntzQuillen} that $\Omega A = T_A(\Omega^1 A)$ the tensor algebra of the $A$-bimodule $\Omega^1 A$. Remark that the standard assumption of \cite{CuntzQuillen} is that $\ell = \C$ the field of complex numbers. However, with minor modifications most results remain valid over an arbitrary basefield and we will refer to statements in \cite{CuntzQuillen} whenever the argument can be repeated verbatim.

\begin{definition}
 A finitely generated $\ell$-algebra $A$ is said to be an {\em $\ell$-qurve} (or {\em quasi-free} \cite{CuntzQuillen} or {\em formally smooth} \cite{KontRos} ) if either of the following two equivalent conditions is satisfied
 \begin{itemize}
 \item{The universal bimodule $\Omega^1_{\ell}(A)$ of derivations is a projective $A$-bimodule.}
 \item{$A$ satisfies the lifting property modulo nilpotent ideals in $\ell-alg$, the category of $\ell$-algebras.}
 \end{itemize}
 Whereas the lifting property extends Grothendieck's characterization of commutative regular algebras (see for example \cite{Iversen}) to the non-commutative setting, such algebras are known to be {\em hereditary} by \cite[Prop. 5.1]{CuntzQuillen} and hence they behave {\em qu}ite like c{\em urves}.
 \end{definition}

 Recall that a finite dimensional $\ell$-algebra $S$ is said to be {\em separable} if and only if $S$ is the direct sum of simple algebras each of which has a center which is a separable field extension of $\ell$. For example, the group algebra $\ell G$ of a finite group $G$ is separable if and only if the order of $G$ is a unit in $\ell$. Separable $\ell$-algebras are known to be $\ell$-qurves by \cite[\S 4]{CuntzQuillen} but should be thought of as corresponding to {\em points}. In fact, they are characterized by either of the following two equivalent conditions
 \begin{itemize}
 \item{$A$ is a projective $A$-bimodule.}
 \item{$A$ satisfies the {\em conjugate} lifting property modulo nilpotent ideals in $\ell-alg$.}
 \end{itemize}
 That is, if $I \triangleleft B$ is a nilpotent ideal and if $\overline{\phi},\overline{\psi} : S \pile{\rTo \\ \rTo}B/I$ are two $\ell$-algebra morphisms which are conjugated by a unit $\overline{b} \in B/I$ then there exist algebra lifts $\phi,\psi : S \pile{\rTo \\ \rTo} B$ and a unit $b \in B$ (mapping to $\overline{b}$) conjugating $\phi$ and $\psi$, see \cite[Prop. 6.1.2]{CuntzQuillen}.
 
 Genuine examples of $\ell$-qurves are : the free algebra $\ell \langle x_1,\hdots,x_m \rangle$, the path algebra $\ell Q$ of a finite quiver $Q$ and the coordinate ring $\ell[C]$ of a smooth affine commutative curve $C$. From these more complicated examples are construed by universal constructions such as taking algebra free products $A \ast A'$ or universal localizations $A_{\Sigma}$. In the next section we will introduce a new class of $\ell$-qurve examples.
 
 For an $\ell$-algebra $A$ recall that the {\em representation scheme} $rep_n~A$ is the affine $\ell$-scheme representing the functor
 \[
 \ell-commalg \rTo sets \qquad \text{defined by} \qquad C \mapsto Hom_{\ell-alg}(A,M_n(C)) \]
 where $\ell-commalg$ is the category of commutative $\ell$-algebras.
 A major motivation for studying $\ell$-qurves comes from the result mentioned in \cite{Kontsevich}, \cite{KontRos} and proved in \cite[(2.2)]{LBnagatn}.
 
 \begin{proposition} \label{smoothness} If $A$ is a $\ell$-qurve, then all representation schemes $rep_n~A$ are smooth affine varieties (possibly having several connected components).
 \end{proposition}
 
 \section{Qurves from graphs}
 
 In this section we will imitate the Bass-Serre theory of the fundamental group of a graph of groups, see \cite{Serre} or \cite{Dicks},  to construct a large class of examples of $\ell$-qurves.
 
 \begin{definition} Let $G = (V,E)$ be a finite graph with vertex-set $V$ and edges $E$. A {\em $G$-graph of $\ell$-qurves} $\mathcal{Q}_G$ is the assignment of 
 \begin{itemize}
 \item{An $\ell$-qurve  $A_v$ to every vertex $v \in V$.}
 \item{A separable $\ell$-algebra $S_e$ to every edge $e \in E$.}
 \item{Inclusions of $\ell$-algebras 
  \[
 S_e \rInto^{i_{e,v}} A_v \quad \text{and} \quad S_e \rInto^{i_{e,w}} A_w \quad \text{for every edge} \quad \xymatrix{\vtx{v} \ar@{-}[r]^e & \vtx{w} }
 \]}
 \end{itemize}
 If, moreover, all vertex-algebras are separable algebras $S_v$ we call this data a {\em $G$-graph of separable algebras} and denote it by $\mathcal{S}_G$.
 \end{definition}
 
In order to construct the {\em fundamental algebra} $\pi_1(\mathcal{Q}_G)$ of a $G$-graph of qurves $\mathcal{Q}_G$ we need to have $\ell$-algebra equivalents for the notions of {\em amalgamated group products} \cite[\S 1.2]{Serre} and of the {\em HNN construction} \cite[\S 1.4]{Serre}. If $S$ is a separable $\ell$-algebra and if $A$ and $A'$ are $S$-algebras, then the {\em coproduct} $A \ast_S A'$ is the algebra representing the functor
 \[
 Hom_{S-\wis{alg}}(A,-) \times Hom_{S-\wis{alg}}(A',-) \]
 in the category $S-alg$ of $S$-algebras, see for example \cite[Chp. 2]{Schofield} for its construction and properties. As for the HNN-construction, let $\alpha,\beta : S \rInto A$ be two $\ell$-algebra embeddings of $S$ in $A$, consider the algebra
 \[
 A \ast_S^{\alpha,\beta} = \frac{A \ast \ell[t,t^{-1}]}{(\beta(s) - t^{-1}\alpha(s)t~:~\forall s \in S)} \]
 
 \begin{lemma} \label{tools} Let $S$ be a separable $\ell$-algebra, $A$ and $A'$ $\ell$-qurves  and $\ell$-embeddings $\alpha,\beta : S \rInto A$ and $S \rInto A'$. Then, the $\ell$-algebras
 \[
 A \ast_S A' \qquad \text{and} \qquad A \ast_S^{\alpha,\beta} \]
 are again $\ell$-qurves.
 \end{lemma}
 
 \begin{proof} Our edge-algebras need to be separable $\ell$-algebras because we will need the conjugate lifting property modulo nilpotent ideals.
 
 A morphism $A \ast_S A' \rTo^g B/I$ is fully determined by morphisms $A \rTo^f B/I$ and $A' \rTo^{f'} B/I$ such that $f | S = f' | S$. As $A$ and $A'$ are quasi-free one has $\ell$-algebra lifts $\tilde{f}  : A \rTo B$ and $\tilde{f'} : A' \rTo B$ whence two morphisms on $S$ which have to be conjugated by an $b \in B^*$ such that $\overline{b} = 1_{B/I}$, that is $f'(s) = b^{-1} f(s) b$ for all $s \in S$. But then, we have a lift $A \ast_S A' \rTo B$ determined by the morphisms $b^{-1} f b$ and $f'$.
 
 A morphism $A \ast_S^{\alpha,\beta} \rTo^g B/I$ determines (and is determined by) a morphism $A \rTo^f B/I$ and a unit $\overline{b} = g(t)$ such that $f \circ \alpha$ and $f \circ \beta : S \rTo B/I$ are conjugated via $\overline{b}$. Because $A$ is quasi-free we have a lift $\tilde{f} : A \rTo B$ and algebra maps $\tilde{f }\circ \alpha$ and $\tilde{f} \circ \beta : S \rTo B$ which reduce to $\overline{b}$ conjugate morphisms. But then there is a unit $b \in B^*$ conjugating $\tilde{f} \circ \alpha$ to $\tilde{f} \circ \beta$ and mapping $t$ to $b$ produces the required lift $A \ast_S^{\alpha,\beta} \rTo B$.
 \end{proof}
 
 However, as often with universal constructions, we have to take care not to end up with the trivial algebra! Because $S$ is semi-simple and $A$ and $A'$ are faithful $S$-algebras it follows from \cite[Chp. 2]{Schofield} that there are inclusions $A \rInto A \ast_S A'$ and $A' \rInto A \ast_S A'$. To prove that $A \rInto A \ast_S^{\alpha,\beta}$ we give another description of the HNN-construction mimicking \cite[\S 1.4]{Serre}. For any $n \in \Z$ take $A[n] \simeq A$ and construct the following amalgamated products
 \[
 A_0 = A, \quad A_1 = A[-1] \ast_S A_0 \ast_S A[1] ,~\hdots~ A_k = A[-k] \ast_S A_{k-1} \ast_S A[k] \]
 with respect to the following embeddings
 \[
 \xymatrix{
A[-1] & & A[0] & & A[1] & & A[2] \\
& S \ar[lu]_{\beta} \ar[ru]^{\alpha}&  & S \ar[lu]_{\beta} \ar[ru]^{\alpha} & & S \ar[lu]_{\beta} \ar[ru]^{\alpha} &}
 \]
 As $S$ is semi-simple we have by \cite[Chp. 2]{Schofield} embeddings $A_0 \subset A_1 \subset A_2 \subset ... $ and hence $A$ embeds in the limit $\tilde{A} = \underset{\rightarrow}{lim} A_n$. The shift-identity 
 \[
 \hdots \rTo A[k-1] \rTo^{id} A[k] \rTo^{id} A[k+1] \rTo \hdots \]
induces an automorphism $\phi$ on $\tilde{A}$ and as the  two algebras below have the same universal property they are isomorphic
 \[
 A \ast_S^{\alpha,\beta} \simeq \tilde{A} [t,t^{-1},\phi] \qquad \text{whence}\qquad A \rInto A \ast_S^{\alpha,\beta} \]
 
 \begin{definition} Let $\mathcal{Q}_G$ be a graph of $\ell$-qurves and let $T$ be a maximal subtree of $G$. We construct the $\ell$-algebra $A_T$ by induction on the number $t$ of edges in $T$. If $t=0$ so $V = \{ v \}$ then $A_T = A_v$.
 If $t > 0$, consider a leaf vertex $v$ with connecting edge $\xymatrix{\vtx{v} \ar@{-}[r]^e & \vtx{w}}$ in $T$. Construct a new tree $T'$ on $t-1$ edges by dropping the vertex $v$ and edge $e$ and construct a new graph of $\ell$-qurves $\mathcal{Q'}_{T'}$ by 
 \[
 A'_w = A_v \ast_{A_e} A_w, \quad A'_{u} = A_u~\quad \text{for $v \not= u \in V$}, \quad A'_f = A_f~\text{for $e \not= f \in E$} \]
 then $A_T \simeq A_{T'}$. Observe that there are embeddings $S_u \rInto^{i_u} A_T$ for every $u \in V$.
 
 Let $G-T = \{ e_1,\hdots,e_r \}$ and take $A_0 \simeq A_T$.
 For every edge $\xymatrix{\vtx{v} \ar@{-}[r]^{e_i} & \vtx{w}}$ in $G-T$ there are two embeddings
 \[
 \alpha_i~:~S_e \rInto^{i_{e_i,v}} S_v \rInto^{i_v} A_{i-1} \quad \text{and} \quad
 \beta_i~:~S_e \rInto^{i_{e_i,w}} S_w \rInto^{i_w} A_{i-1} \]
and we define
 \[
 A_i \simeq A_{i-1} \ast_{S_e}^{\alpha_i,\beta_i} \]
The algebra $A_r$ is then called the {\em fundamental algebra of the graph of $\ell$-qurves $\mathcal{Q}_G$} and is denoted by $\pi_1(\mathcal{Q}_G)$.
 \end{definition}
 
 \begin{theorem} If $\mathcal{Q}_G$ is a graph of $\ell$-qurves, the fundamental algebra $\pi_1(\mathcal{Q}_G)$ is again an $\ell$-qurve.
 \end{theorem}
 
 \begin{proof} Immediate from the construction and lemma~\ref{tools}.
 \end{proof}

 \section{Qurve group algebras}
 
 The classification of $\ell$-qurves is way out of reach at the moment so it is important to have partial classifications. In \cite[\S 6]{CuntzQuillen} the finite dimensional $\ell$-qurves were shown to be the hereditary finite dimensional $\ell$-algebras (and hence Morita equivalent to path algebras $\ell Q$ of a finite quiver $Q$ without oriented cycles). In this section we will classify the group algebras $\ell H$ for $H$ a finitely generated group which are $\ell$-qurves. The desired answer is that these are precisely the $\ell H$ with $H$ a {\em virtually free group} (that is, $H$ has a free subgroup of finite index) but we have to take the characteristic of $\ell$ into account (observe that finite groups are virtually free).
 
 If $\mathcal{G}_G$ is a graph of {\em finite groups} as in \cite{Serre} such that all orders are invertible in $\ell$, then we can associate to it a graph of separable $\ell$-algebras $\mathcal{S}_G$ by taking
 \[
 S_v = \ell G_v~\quad \forall v \in V \qquad \text{and} \qquad S_e = \ell G_e~\quad \forall e \in E \]
 with embeddings determined by the group-embeddings. If $\pi_1(\mathcal{G}_G)$ is the {\em fundamental group} of $\mathcal{G}_G$ as in \cite[\S 5.1]{Serre} then the point of the construction in the previous section is that
 \[
 \ell \pi_1(\mathcal{G}_G) \simeq \pi_1(\mathcal{S}_G) \]
 and hence these group algebras are $\ell$-qurves. The connection with virtually free groups is provided by a result of Karrass, see for example \cite[Thm. 3.5]{Wall}. The following statements are equivalent for a finitely generated group $H$
 \begin{itemize}
 \item{$H = \pi_1(\mathcal{G}_G)$ for a graph of finite groups.}
 \item{$H$ is a virtually free group.}
 \end{itemize}
 For example, all congruence subgroups in the modular group $SL_2(\Z)$ are virtually free. On the other hand, the third braid group $B_3 = \langle s,t~|~s^2=t^3 \rangle$ is not virtually free. Note that very little is known about simple representations of congruence subgroups. For some low dimensional classifications of $SL_2(\Z)$-representations see \cite{Wenz}. 
 
 \begin{theorem} The following statements are equivalent for a finitely generated group $H$:
 \begin{enumerate}
 \item{The group algebra $\ell H$ is an $\ell$-qurve.}
 \item{$H$ is a virtually free group such that in a description $H = \pi_1(\mathcal{G}_G)$ all orders of the vertex groups $G_v$ are finite and invertible in $\ell$.}
 \end{enumerate}
 \end{theorem}
 
 \begin{proof} If $\ell H$ is a quasi-free $\ell$-algebra, it has to be hereditary by \cite[Prop. 5.1]{CuntzQuillen} and hence, in particular, its augmentation ideal $\omega_H$ mast be a projective left $\ell H$-module. By a result of Dunwoody, see \cite[Thm. IV.2.12]{Dicks} this is equivalent to $H$ being the fundamental group of a graph of finite groups $\mathcal{G}_G$ such that all vertex-group orders are invertible in $\ell$, whence (2) follows. The converse implication follows from the discussion preceding the statement and the last section.
 \end{proof}
 
 If $char(\ell) = 0$ it follows from this and proposition~\ref{smoothness} that all representation schemes $rep_n~\ell H$ are smooth affine varieties whenever $H$ is a finitely generated virtually free group.

 \section{The component semigroup}
 
 From now on we will assume that $\ell = \overline{\ell}$ is algebraically closed. In the appendix we will replace the component semigroup by a component co-algebra over an arbitrary basefield $\ell$. If $A$ is an $\overline{\ell}$-qurve we know from proposition~\ref{smoothness} that all representation schemes are smooth affine varieties.
 
 \begin{definition} For an $\overline{\ell}$-qurve $A$ the smooth variety $rep_n~A$ decomposes into connected (equivalently, irreducible) components
 \[
 rep_n~A = \bigsqcup_{| \alpha | = n} rep_{\alpha}~A \]
 where $\alpha$ is a label. We call $\alpha$ a {\em dimension vector} of {\em total dimension} $| \alpha | = n$.
 \end{definition}
 
 An $\overline{\ell}$-point of $rep_n~A$ is an $n$-dimensional left $A$-module and the direct sum of modules defines the {\em sum maps}
 \[
 rep_n~A \times rep_m~A \rTo rep_{n+m}~A \]
 If we decompose these varieties into their connected components and use the fact that the image of two connected varieties is again connected, we can define a semigroup.
 
 \begin{definition} The {\em component semigroup} $comp(A)$ is the set of all dimension vectors $\alpha$ equipped with the addition $\alpha +\beta = \gamma$ where $\gamma$ determines the unique component $rep_{\gamma}~A$ of $rep_{n+m}~A$ containing the image of $rep_{\alpha}~A \times rep_{\beta}~A$ under the sum map
 \[
 \bigsqcup_{| \alpha | = n} rep_{\alpha}~A \times \bigsqcup_{| \beta | = m} rep_{\beta}~A \rTo \bigsqcup_{| \gamma| = n+m} rep_{\gamma}~A \]
 $comp(A)$ is a commutative semigroup with an augmentation map $comp(A) \rTo \N$ sending a dimension vector $\alpha$ to its total dimension $| \alpha |$.
 \end{definition}
 
 \noindent
 Here are some examples :
\begin{itemize}
\item{For $A = M_{n_1}(\overline{\ell}) \oplus \hdots \oplus M_{n_k}(\overline{\ell})$ semi-simple, $comp(A) = (\N n_1,\hdots,\N n_k) \subset \N^k$.}
\item{For $A = \overline{\ell} Q$ a path algebra we have $comp(A) = \N^k$ where $k$ is the number of vertices of the quiver $Q$.}
\item{For a direct sum $A = A_1 \oplus A_2$ we have $comp(A) = comp(A_1) \oplus comp(A_2)$.}
\item{For a free algebra product $A = A_1 \ast A_2$ we have that $comp(A_1)$ is the fibered product (using the augmentation) $comp(A_1) \times_{\N} comp(A_2)$, see \cite[Prop. 1]{Morrison}.}
\end{itemize}
In \cite[Question 2]{Morrison} K. Morrison asked whether $comp(A)$ is always a free Abelian semigroup (as in the examples above). However, even for $A$ an $\overline{\ell}$-qurve, reality is more complex as one can remove components by the process of universal localization (see for example \cite{Schofield} for definition and properties of universal localization).

\begin{proposition} For every sub semigroup $S \subset \N$, there is an $\overline{\ell}$-qurve $A$ with
\[
comp(A) = S \]
as augmented semigroups.
\end{proposition}

\begin{proof} Suppose first that $gcd(S)=1$, that is the elements of $S$ are coprime. By using results on polynomial- and rational identities of matrices (see for example \cite{Rowen}) it was proved in \cite{LBbergman} that there is an affine $\ell$-algebra with presentation
\[
A = \frac{\overline{\ell} \langle x_1,\hdots,x_a,y_1,\hdots,y_b \rangle}{(1-y_i p_i(x_1,\hdots,x_a,y_1,\hdots,y_{i-1})~:~1 \leq i \leq b)}
\]
(with each of the $p_i \in \ell \langle x_1,\hdots,x_a,y_1,\hdots,y_{i-1} \rangle$) having the property that $A$ has finite dimensional representations of dimensions exactly the elements of $S$. $A$ is a universal localization of $\overline{\ell} \langle x_1,\hdots,x_a \rangle$ and hence is an $\overline{\ell}$-qurve (for example use \cite[Thm. 10.6]{Schofield} to prove that $\Omega^1(A)$ is a projective $A$-bimodule). As such, for every $n$, $rep_n~A$ is a Zariski open subset (possibly empty) of $rep_n~\overline{\ell} \langle x_1,\hdots,x_a \rangle = M_n(\overline{\ell})^{\times a}$ and is therefore irreducible (when non-empty). Therefore,
$comp(A) = S \subset \N$ and consists precisely of those $n \in \N$ for which none of the $p_i$ (when expressed as a rational non-commutative function in $x_1,\hdots,x_a$) is a rational identity for $n \times n$ matrices.

For the general case, assume that $gcd(S) = m$ and take $S' = S/m$ with associated algebra (as above) $A'$ for which $comp(A') = S' \subset \N$. But then,
\[
comp(A' \ast M_m(\overline{\ell})) = S' \times_{\N} \N m = S \]
and $A = A' \ast M_m(\overline{\ell})$ is again an $\overline{\ell}$-qurve.
\end{proof}

\section{Tits and Euler forms}

In this section we will define bilinear forms on $comp(A)$ (when $A$ is an $\overline{\ell}$-qurve) generalizing the Tits- and Euler-forms on the dimension vectors of a quiver. Let $rep~A$ be the Abelian category of all finite dimensional representations of $A$. If $A$ is an affine $\overline{\ell}$-algebra, then
$Hom_A(M,N)$ and $Ext^1_A(M,N)$ are finite dimensional $\overline{\ell}$-spaces for all $M,N \in rep~A$.

If $A$ is hereditary (for example, if $A$ is an $\overline{\ell}$-qurve) we have that $\chi_A(M,-)$ and $\chi_A(-,N)$ are additive on short exact sequences in $rep~A$ where
\[
\chi_A(M,N) = dim_{\overline{\ell}} Hom_A(M,N) - dim_{\overline{\ell}} Ext^1_A(M,N) \]
For $M \in rep~A$ define its {\em semi-simplification} $M^{ss}$ to be the semi-simple $A$-module obtained by taking the direct sum of the Jordan-H\"older components of $M$. From additivity on short exact sequences it follows for all $M,N \in rep~A$ that
\[
\chi_A(M,N) = \chi_A(M^{ss},N^{ss}) \]
For $\alpha,\beta \in comp(A)$ it follows from \cite{KraftLNM} and \cite[lemma 4.3]{CrawleySchroer} that the functions
\[
rep_{\alpha}~A \times rep_{\beta}~A \rTo \Z \qquad (M,N) \mapsto \begin{cases}
dim_{\overline{\ell}} Hom_A(M,N) \\
dim_{\overline{\ell}} Ext^1_A(M,N)
\end{cases}
\]
are upper semicontinuous. In particular, there are Zariski open subsets (whence dense by irreducibility)
of $rep_{\alpha}~A \times rep_{\beta}~A$ where these functions attain a minimum. Following \cite{Schofieldgeneric} we will denote these minimal values by $hom(\alpha,\beta)$ resp. $ext(\alpha,\beta)$. 

The group $GL_n$ acts on $rep_n~A$ by base-change and orbits $\Oscr(M)$ under this action are precisely the isomorphism classes of $n$-dimensional left $A$-modules. From \cite{Gabriel} we recall that the semi-simplification $M^{ss}$ belongs to the Zariski closure $\overline{\Oscr(M)}$ of the orbit and that $Ext^1_A(M,M)$ can be identified to the {\em normal space} to the orbit $\Oscr(M)$ with respect to the scheme structure on $rep_n~A$.

\begin{proposition} \label{tits} Let $A$ be an affine $\overline{\ell}$-algebra.
\begin{enumerate}
\item{If $rep_{\gamma}~A$ is a smooth variety, then for all $M \in rep_{\gamma}~A$ we have
\[
|\gamma|^2 - \chi_A(M,M) = dim~rep_{\gamma}~A \]
and hence $\chi_A(M,M)$ is constant on $rep_{\gamma}~A$.}
\item{If $rep_{\alpha}~A$, $rep_{\beta}~A$ and $rep_{\alpha+\beta}~A$ are smooth varieties, then
\[
\chi_A(M,N) + \chi_A(N,M) \]
is a constant function on $rep_{\alpha}~A \times rep_{\beta}~A$.}
\end{enumerate}
\end{proposition}

\begin{proof}
If $rep_{\gamma}~A$ is smooth in $M$, it follows from the above remarks that
\[
T_M rep_{\gamma}~A = Ext^1_A(M,N) \oplus T_M \Oscr(M), \qquad \Oscr(M) = GL_{| \gamma |} / Stab(M) \]
where $Stab(M)$ is the stabilizer subgroup which by \cite{KraftLNM} has the same dimension as $Hom_A(M,M)$. Therefore,
\[
\begin{split}
dim~rep_{\gamma}~A &= dim_{\overline{\ell}} T_M rep_{\gamma}~A \\ &= dim_{\overline{\ell}} Ext^1_A(M,M) + | \gamma |^2 - dim_{\overline{\ell}} Hom_A(M,M)
\end{split}
 \]
whence (1). (2) follows from this by considering the point $M \oplus N \in rep_{\alpha+\beta}~A$ and using bi-additivity of $\chi_A$.
\end{proof}

\begin{definition} If $A$ is an $\overline{\ell}$-qurve, then for all $\alpha \in comp(A)$ the representation variety $rep_{\alpha}~A$ is smooth. Therefore, the constant value 
\[
(\alpha,\beta)_A = \chi_A(M,N) + \chi_A(N,M) \]
on $rep_{\alpha}~A \times rep_{\beta}~A$ defines a symmetric bilinear form
\[
(-,-)_A~:~comp(A) \times comp(A) \rTo \Z \]
which we call the {\em Tits-form} of the $\overline{\ell}$-qurve $A$.
\end{definition}

\noindent
For general affine $\overline{\ell}$-algebras $\chi_A(M,N)+\chi_A(N,M)$ does not have to be constant and the foregoing result can be used to deduce singularity of specific representation varieties.

\begin{example} Let $A = \overline{\ell} B_3$ be the group-algebra of the {\em third braid group}
$B_3 = \langle s,t~|~s^2=t^3 \rangle$.
The one dimensional representation variety is the cusp minus the singular origin
\[
rep_1~A = \{ (x,y) \in \overline{ \ell}^2~|~x^3=y^2 \} - \{ (0,0) \} \]
and hence is a smooth affine variety. As all points are simple $A$-modules we have that $dim_{\overline{\ell}}~Hom_A(-,-)$ is equal to zero on the open set $rep_1~A \times rep_1~A - \Delta$ and is equal to one on the {\em diagonal} $\Delta$. As for $dim_{\overline{\ell}}~Ext^1_A(-,-)$ this is zero on
$rep_1~A \times rep_1~A -( \Delta \sqcup \Delta_1 \sqcup \Delta_2)$ where
\[
\begin{cases}
\Delta_1 &= \{ ((x,y),(\rho x,-y))~:~x^3=y^2 \} \\ \Delta_2 &= \{ ((x,y),(\rho^2 x,-y))~:~x^3=y^2 \}
\end{cases}  \]
for $\rho$ a primitive third root of unity. As a consequence, $\chi_A(M,N)$ is zero on the Zariski open subset $rep_1~A \times rep_1~A - (\Delta_1 \sqcup \Delta_2)$ and is equal to $-1$ on $\Delta_1 \sqcup \Delta_2$. Therefore, $\overline{\ell} B_3$ is not an $\overline{\ell}$-qurve. In fact, $rep_2~\overline{\ell} B_3$ is not smooth.
\end{example}

If $\alpha$ is the dimension vector of a simple representation of $A$, then there is a  Zariski open subset $simp_{\alpha}~A$ of simple representations in $rep_{\alpha}~A$.

\begin{proposition} If $A$ is an $\overline{\ell}$-qurve and $\alpha,\beta$ are dimension vectors of simple representations, then the function
\[
\chi_A(S,T) \]
is constant on $simp_{\alpha}~A \times simp_{\beta}~A$.
\end{proposition}

\begin{proof} There is a Zariski open subset $U \subset simp_{\alpha}~A \times simp_{\beta}~A$ consisting of couples $(S',T')$ such that
\[
dim_{\overline{\ell}} Ext^1_A(S',T') = ext(\alpha,\beta) \quad \text{and} \quad
dim_{\overline{\ell}} Ext^1_A(T',S') = ext(\beta,\alpha) \]
Hence, for all $(S,T) \in simp_{\alpha}~A \times simp_{\beta}~A$
\[
\begin{cases}
dim_{\overline{\ell}} Ext^1_A(S,T) \geq dim_{\overline{\ell}} Ext^1_A(S',T') \\
dim_{\overline{\ell}} Ext^1_A(T,S) \geq dim_{\overline{\ell}} Ext^1_A(T',S')
\end{cases}
\]
If $\alpha \not= \beta$ (or if $\alpha=\beta$ and $S \not\simeq T$) $\chi_A(S,T) = - dim_{\overline{\ell}} Ext^1_A(S,T)$ and hence the above inequalities must be equalities by proposition~\ref{tits}. Remains to prove for $S,T \in simp_{\alpha}~A$ with $S \not\simeq T$ that $\chi_A(S,S) = \chi_A(S,T)$. Consider the two semi-simple representations $M=S \oplus S$ and $N = S \oplus T$ in $rep_{2 \alpha}~A$. From proposition~\ref{tits}~(1) we get
\[
\begin{split}
4 \chi_A(S,S) &= \chi_A(S,S) + \chi_A(T,T) + \chi_A(S,T) + \chi_A(T,S) \\
&= 2 \chi_A(S,S) + 2 \chi_A(S,T) \\
\end{split}
\]
(using proposition~\ref{tits}~(1) and the above fact that $\chi_A(S,T) = \chi_A(T,S)$) whence $\chi_A(S,S) = \chi_A(S,T)$.
\end{proof}

If $M \in rep~A$, its semi-simplification has as isotypical decomposition
\[
M = S_1^{\oplus e_1} \oplus \hdots \oplus S_k^{\oplus e_k} \]
with all $S_i$ non-isomorphic. If $S_i \in rep_{\beta_i}~A$ we say that the {\em representation type} of $M$ (which is determined upto permutation of the $(e_i,\beta_i)$ terms).
\[
\tau(M) = (e_1,\beta_1;\hdots;e_k,\beta_k) \]

\begin{proposition} \label{euler} If $A$ is an $\overline{\ell}$-qurve, the {\em Euler-form}
\[
\chi_A(M,N) = dim_{\overline{\ell}} Hom_A(M,N) - dim_{\overline{\ell}} Ext^1_A(M,N) \]
depends only on the representation types $\tau(M)$ and $\tau(N)$.
\end{proposition}

\begin{proof}
Follows from the foregoing result by observing that $\chi_A(M,N) = \chi_A(M^{ss},N^{ss})$.
\end{proof}

In particular, there is a Zariski open subset in $rep_{\alpha}~A \times rep_{\beta}~A$ of couples $(M,N)$ on which the value of $\chi_A(M,N)$ is constant and equal to the {\em Euler form}
\[
\chi_A(\alpha,\beta) = hom(\alpha,\beta) - ext(\alpha,\beta) \]
Clearly, this open set contains all representations of {\em generic representation type} $\tau_{gen}$, see for example \cite{LBProcesi}. In fact, if $char(\overline{\ell}) = 0$ the proof of proposition~\ref{exts} implies that $\chi_A(M,N)$ is constant on $rep_{\alpha}~A \times rep_{\beta}~A$.

\section{One quiver to rule them all}

If $A$ is an $\overline{\ell}$-qurve, we will denote with $\Sigma_A$ the minimal set of semigroup-generators of the component semigroup $comp(A)$. Observe that $\Sigma_A$ is well-defined as it follows from the Jordan-H\"older decomposition that
\[
\Sigma_A = \{ \alpha \in comp(A)~|~simp_{\alpha}~A = rep_{\alpha}~A \} \]
In particular, it follows from proposition~\ref{euler} that $\chi_A(S,T) = \chi_S(\alpha,\beta)$ for all representations $S \in rep_{\alpha}~A$ and $T \in rep_{\beta}~A$ if $\alpha,\beta \in \Sigma_A$. In all examples known to us, $\Sigma_A$ is a finite set.

\begin{definition} If $A$ is an $\overline{\ell}$-qurve, we define its {\em one-quiver} $Q_1(A)$ to be the quiver on the (possibly infinite) vertex set $\{ v_{\alpha}~|~\alpha \in \Sigma_A \}$ such that the number of directed arrows (loops) from $v_{\alpha}$ to $v_{\beta}$ is given by
\[
\#~\{~\xymatrix{\vtx{\alpha} \ar[r] & \vtx{\beta}}~\} = \delta_{\alpha \beta} - \chi_A(\alpha,\beta) \]
The {\em one-dimension vector} $\alpha_1(A)$ for $A$ is the dimension vector for $Q_1(A)$ having as its $v_{\alpha}$-component the total dimension $| \alpha |$.
\end{definition}

If $Q_1(A)$ is a quiver on finitely many vertices $\{ v_1,\hdots,v_k \}$ and $\alpha_1(A) = (n_1,\hdots,n_k)$, we can define the $\overline{\ell}$-algebra
\[
B(Q_1(A),\alpha_1(A)) = \begin{bmatrix}
B_{11} & \hdots & B_{1k} \\
\vdots & & \vdots \\
B_{k1} & \hdots & B_{kk} \end{bmatrix}
\]
where $B_{ij}$ is the $n_i \times n_j$ block matrix having all its components equal to the sub vectorspace of the path algebra $\overline{\ell} Q_1(A)$ spanned by all oriented paths in $Q_1(A)$ starting at vertex $v_i$ and ending in $v_j$. Observe, that $B(Q_1(A),\alpha_1(A))$ is Morita equivalent to the path algebra $\overline{\ell} Q_1(A)$ and as such is again an $\overline{\ell}$-qurve.

\begin{example}[Deligne-Mumford curves]  Recall from \cite[Coroll. 7.8]{ChanIngalls} that a {\em smooth Deligne-Mumford curve} which is generically a scheme, determines (and is determine by) a smooth affine curve $X$ and an {\em hereditary order} $A$ over $\overline{\ell} [X]$. As such, $A$ is an $\overline{\ell}$-qurve with center $\overline{\ell} [X]$ and is a subalgebra of $M_n(\overline{\ell} (X))$ for some $n$  called the p.i.-degree of $A$. If $\mathfrak{m}_x$ is the maximal ideal of $\overline{\ell}[X]$ corresponding to the point $x \in X$ then for all but finitely many exceptions $\{ x_1,\hdots,x_l \}$ we have that
\[
A / \mathfrak{m}_x A \simeq M_n(\overline{\ell)} \]
For the exceptional points (the {\em ramification locus} of $A$) there are finitely many maximal ideals $\{ P_1(i),\hdots,P_{k_i}(i) \}$ of $A$ lying over $\mathfrak{m}_{x_i}$ and
\[
A/P_j(i) \simeq M_{n_j(i)}(\overline{\ell}) \qquad \text{with} \qquad n_1(i) + \hdots + n_{k_i}(i) = n \]
As a consequence, $rep_l~A$ for all $l < n$ consists of finitely many closed orbits each corresponding to a maximal ideal $P_j(i)$ such that $A/P_j(i) \simeq M_l(\overline{\ell})$. Hence, the component semigroup $comp(A)$ has generators $\alpha_j(i)$ for all $1 \leq i \leq l$ and $1 \leq j \leq k_i$ and relations for all $1 \leq i,j \leq l$
\[
\alpha_1(i)+\hdots+\alpha_{k_i}(i) = \alpha_1(j) + \hdots + \alpha_{k_j}(j) \]
From direct calculation or using \cite[Prop. 6.1]{LBlocalstructure} it follows that the one quiver $Q_1(A)$ is the disjoint union of $l$ quivers of type $\tilde{A}_{k_i}$, that is the $i$-th component is $Q_1(A)(i)$ and is the quiver on $k_i$ vertices
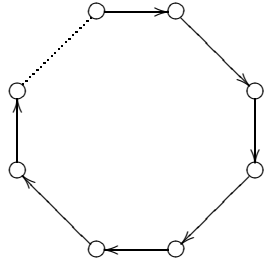
\[
\xymatrix{
& \vtx{} \ar[r] & \vtx{} \ar[dr] & \\
\vtx{}\ar@{.}[ur] & &  & \vtx{} \ar[d]  \\
\vtx{}\ar[u] & &  & \vtx{}\ar[dl] \\
& \vtx{} \ar[ul] & \vtx{} \ar[l] & }
\]
and the corresponding components for the one dimension vector $\alpha_1(A)$ are $\alpha_1(A)(i) = (n_1(i),\hdots,n_{k_i}(i))$. Therefore, the associated algebra
\[
B(Q_1(A),\alpha_1(A)) = B_1 \oplus \hdots \oplus B_l \]
where $B_i$ is the block-matrix algebra
\[
\begin{bmatrix}
M_{n_1(i) \times n_1(i)}(\overline{\ell} [x]) & M_{n_1(i) \times n_2(i)}(\overline{\ell} [x]) & \hdots & M_{n_1(i) \times n_{k_i}(i)}(\overline{\ell} [x]) \\
M_{n_2(i) \times n_1(i)}( x \overline{\ell} [x]) & M_{n_2(i) \times n_2(i)}(\overline{\ell} [x]) & \hdots & M_{n_2(i) \times n_{k_i}(i)}(\overline{\ell} [x]) \\
\vdots & \vdots & & \vdots \\
M_{n_{k_i}(i) \times n_1(i)}(x \overline{\ell} [x]) & M_{n_{k_i}(i) \times n_2(i)}(x \overline{\ell} [x]) & \hdots & M_{n_{k_i}(i) \times n_{k_i}(i)}(\overline{\ell} [x]) 
\end{bmatrix}
\]
It follows from \cite[Chp. 9]{Reiner} or \cite[Prop. 6.1]{LBlocalstructure} that in a neighborhood of $x_i$ the $\overline{\ell}$-qurve $A$ is \'etale isomorphic to $B_i$.
\end{example}

Elsewhere, we will generalize this example by relating the $\overline{\ell}$-qurve $A$ with the algebra $B(Q_1(A),\alpha_1(A))$ using the formal tubular neighborhood theorem \cite[\S 6]{CuntzQuillen}. Here, we will use the {\em one-quiver-setting} $(Q_1(A),\alpha_1(A))$ to describe the $GL_n$-\'etale local structure of $rep_n~A$ in the neighborhood of a semi-simple representation. As this description uses the Luna slice result, we will assume that $char(\overline{\ell}) = 0$ in the remainder of this section. We recall the construction of the {\em local quiver} and refer to \cite{LBnagatn} and \cite{LBlocalstructure} for details and proofs.

\begin{definition}
Let $M \in rep_{\alpha}~A$ be a semi-simple $A$-module of representation type $\tau_M = (e_1,\gamma_1;\hdots;e_l,\beta_l)$, that is
\[
M = S_1^{\oplus e_1} \oplus \hdots \oplus S_l^{\oplus e_l} \]
with all $S_i$ non-isomorphic and of dimension vector $\gamma_i$. 

The {\em local quiver} $Q_M$ is the quiver on $l$ vertices (corresponding to the distinct simple components of $M$) such that the number of directed arrows from $v_i$ to $v_j$ is equal to $dim_{\overline{\ell}}~Ext^1_A(S_i,S_j)$.

The {\em local dimension vector} $\alpha_M  = (e_1,\hdots,e_l)$ determined by the multiplicities $e_i$ of the simple components of $M$. 
\end{definition}

 Observe that we know already that the quiver $Q_M$ only depends on the representation type $\tau_M$ of $M$ and not on the choice of the simple components $S_i$.  The relevance of this {\em local quiver setting} $(Q_M,\alpha_M)$ is that it determines the $GL_n$-equivariant \'etale structure of $rep_{\alpha}~A$ in a neighborhood of the closed orbit $\Oscr(M)$ by the results from \cite{LBnagatn}. 

As $n = \sum_i e_i | \gamma_i |$ there is an embedding of $GL(\alpha_M)$ into $GL_n$ and with respect to this embedding there is a $GL_n$-equivariant \'etale isomorphism between
\begin{itemize}
\item{$rep_{\alpha}~A$ in a neighborhood of $\Oscr(M)$, and}
\item{$GL_n \times^{GL(\alpha_M)} rep_{\alpha_M}~Q_M$ is a neighborhood of $\Oscr(1_n,0)$}
\end{itemize}
 where $0$ is the zero representation. We will show that the one-quiver setting $(Q_1(A),\alpha_1(A))$ contains enough information to describe all these local quiver settings $(Q_M,\alpha_M)$ whenever $A$ is an $\overline{\ell}$-qurve.

$\Sigma_A = \{ \beta_i ~|~i \in I \}$ is the set of semigroup generators of $comp(A)$ (possibly infinite). For any $\alpha \in comp(A)$ we can write
\[
\alpha = a_1 \beta_1 + \hdots + a_k \beta_k \qquad a_i \in \N \]
(possibly in many several ways) with the $\beta_i \in \Sigma_A$. If the set of vertices $V \leftrightarrow \Sigma_A$ is infinite, we can always replace the infinite one-quiver setting $(Q_1(A),\alpha_1(A))$ by a finite quiver setting $(supp(\alpha),\alpha_1(A) | supp(\alpha))$ where $supp(\alpha)$ is the {\em support} of $\alpha$, that is those vertices $\beta_i \in V \leftrightarrow \Sigma_A$ such that $a_i \in \N_+$ in a fixed description of $\alpha$ in terms of the semigroup generators. For notational reasons, we denote this finite quiver setting again by $(Q_1(A),\alpha_1(A))$.

\begin{proposition} \label{simpele} The one-quiver setting $(Q_1(A),\alpha_1(A))$ contains enough information to determine $simp(A)$ the set of all dimension vectors of simple finite dimensional representations of $A$.
\end{proposition}

\begin{proof} 
If $\alpha \in comp(A)$, fix a description
\[
\alpha = a_1 \beta_1 + \hdots + a_k \beta_k \]
with $a_i \in \N_+$ and $\{ \beta_1,\hdots,\beta_k \}$ among the semigroup generators of $comp(A)$. This implies that there are points in $rep_{\alpha}~A$ corresponding to semi-simple representations
\[
M = S_1^{\oplus a_1} \oplus \hdots \oplus S_k^{\oplus a_k} \]
where the $S_i$ are distinct simple representations in $rep_{\beta_i}~A$. But then the local quiver setting of $M$ in $rep_{\alpha}~A$, $(Q_M,\alpha_M)$  is just $(Q_1(A),\epsilon)$ where $\epsilon=(a_1,\hdots,a_k)$. Because $rep_{\alpha}~A$ is irreducible, it follows that $\alpha \in simp(A)$ if and only if $\epsilon$ is the dimension vector of a simple representation of $Q_1(A)$. These dimension vectors have been classified in \cite{LBProcesi} and we recall the result.

Let $\chi$ be the Euler-form of $Q_1(A)$, that is $\chi = (c_{ij})_{i,j} \in M_k(\Z)$ with $c_{ij} = \delta_{ij} - \# \{ \xymatrix{\vtx{i} \ar[r] & \vtx{j} }~\}$ and let $\delta_i$ be the dimension vector of a vertex-simple concentrated in vertex $v_i$. Then, $\epsilon$ is the dimension vector of a simple representation of $Q_A$ if and only if the following conditions are satisfied : 
\begin{enumerate}
\item{the support $supp(\epsilon)$ is a strongly connected subquiver of $Q_A$ (there is an oriented cycle in $supp(\epsilon)$ containing each pair $(i,j)$ such that $\{ v_i,v_j \} \subset supp(\epsilon)$)}
\item{for all $v_i \in supp(\epsilon)$ we have the numerical conditions
\[
\chi(\epsilon,\delta_i) \leq 0 \qquad \text{and} \qquad \chi(\delta_i,\epsilon) \leq 0 \]
{\em unless} $supp(\epsilon)$ in an oriented cycle of type $\tilde{A}_l$ for some $l$ in which case all components of $\epsilon$ must be equal to one.}
\end{enumerate}
The statement follows from this.
\end{proof}

\begin{proposition} \label{exts} The one-quiver setting $(Q_1(A),\alpha_1(A))$ contains enough information to compute the $\overline{\ell}$-dimension of $Ext^1_A(S,T)$ for all finite dimensional simple representations $S$ and $T$ of $A$. 

If $S \in rep_{\alpha}~A$ where $\alpha = \sum_i a_i \beta_i$ and $T \in rep_{\beta}~A$ where $\beta = \sum_i b_i \beta_i$, then
\[
dim_{\overline{\ell}}~Ext^1_A(S,T) = -\chi_{Q_1(A)}(\epsilon,\eta) \]
for $\epsilon = (a_1,\hdots,a_k)$ and $\eta = (b_1,\hdots,b_k)$.
\end{proposition}

\begin{proof}
Let $S_i$ and $T_i$ be distinct simples in $rep_{\beta_i}~A$ and consider the semi-simple representations $M$ resp. $N$ in $rep_{\alpha}~A$ resp. $rep_{\beta}~A$
\[
M = S_1^{\oplus a_1} \oplus \hdots \oplus S_k^{\oplus a_k} \quad \text{and} \quad
N = T_1^{\oplus b_1} \oplus \hdots \oplus T_k^{\oplus b_k} \]
By the foregoing proposition, we have complete information on the local quiver setting of $M \oplus N$
in $rep_{\alpha+\beta}~A$ from $(Q_1(A),\alpha_1(A))$. 
By assumption on $\alpha$ and $\beta$ there is a Zariski open subset of simples $S' \in rep_{\alpha}~A$ and simples $T' \in rep_{\beta}~A$ such that $S' \oplus T'$ lies in a neighborhood of $M \oplus N$. 

It follows from \cite{LBProcesi} that one can reconstruct the local quiver setting of $S' \oplus T'$ from that of $M \oplus N$. This local quiver has two vertices $\{ v_1,v_2 \}$ with $-\chi_Q(\eta,\epsilon)$ arrows from $v_1$ to $v_2$ and $-\chi_Q(\epsilon,\eta)$ arrows from $v_2$ to $v_1$. In $v_1$ there are $1-\chi_Q(\epsilon,\epsilon)$ loops and in $v_2$ there are $1-\chi_Q(\eta,\eta)$ loops. The dimension vector is $(1,1)$.
From this we deduce that
\[
dim_{\overline{\ell}}~Ext^1_A(S',T')  = - \chi(\epsilon,\eta) \]
but we have seen before that  the extension-dimension only depends on the representation type and not on the choice of simples, hence this number is also equal to $dim_{\overline{\ell}}~Ext^1_A(S,T)$.
\end{proof}

\begin{theorem} The one-quiver setting $(Q_1(A),\alpha_1(A))$ contains enough information to construct the local quiver setting $(Q_M,\alpha_M)$ for every semi-simple representation
\[
M = S_1^{\oplus e_1} \oplus \hdots \oplus S_l^{\oplus e_l} \]
of $A$.
\end{theorem}

\begin{proof}
This is a direct consequence of the foregoing two propositions. To begin, we can determine the possible dimension vectors $\alpha_i$ of the simple components $S_i$. Write $\alpha_i = \sum_{j=1}^k a_j(i) \beta_j$ then $\epsilon_i = (a_1(i),\hdots,a_k(i))$ must be the dimension vector of a simple representation of the associated quiver $Q_1(A)$. Moreover, by the previous theorem we know that
\[
dim_{\overline{\ell}}~Ext^1_A(S_i,S_j) = \delta_{ij}-\chi(\epsilon_i,\epsilon_j) \]
and hence we have full knowledge of the local quiver $Q_M$.
\end{proof}

\section{The one-quiver for $\pi_1(\mathcal{S}_G)$}

In this section we will construct the one-quiver setting for the fundamental algebra $\pi_1(\mathcal{S}_G)$ of a graph $\mathcal{S}_G$ of separable (that is, semi-simple) $\overline{\ell}$-algebras. As an intermediary step we will construct a finite quiver $Q_0(\mathcal{S}_G)$ such that finite dimensional representations of $\pi_1(\mathcal{S}_G)$ correspond to certain finite dimensional representations  of the path algebra $\overline{\ell} Q_0(\mathcal{S}_G)$.

 We have decomposition of the vertex- and edge-algebras
 \[
 S_v = M_{d_v(1)}(\overline{\ell}) \oplus \hdots \oplus M_{d_v(n_v)}(\overline{\ell}) \quad \text{resp.} \quad
 S_e = M_{d_e(1)}(\overline{\ell}) \oplus \hdots \oplus M_{d_e(n_e)}(\overline{\ell}) \]
 The embeddings $S_e \rInto S_v$ are depicted via Bratelli-diagrams or, equivalently, by natural numbers $a^{(ev)}_{ij}$ for $1 \leq i \leq n_e$ and $1 \leq j \leq n_v$ satisfying the numerical restrictions
 \[
 d_v(j) = \sum_{i=1}^{n_e} a^{(ev)}_{ij} d_e(i) \qquad \text{for all $1 \leq j \leq n_v$ and all $v \in V$ and $e \in E$} \]
 Remark that these numbers give the {\em restriction data}, that is, the multiplicities of the simple components of $S_e$ occurring in the restriction $V^{(v)}_j \downarrow_{S_e}$ for the simple components $V_j^{(v)}$ of $S_v$. From these decompositions and Schur's lemma it follows that for any edge
 $\xymatrix{\vtx{v} \ar@{-}[r]^e & \vtx{w} }$ in the graph $G$ we have
 \[
 Hom_{S_e}(V_i^{(v)},V_j^{(w)}) = \sum_{k=1}^{n_e} a_{ki}^{(ev)} a_{kj}^{(ew)} = n_{ij}^{(e)} \]
 
 \begin{definition}
 For a graph $\mathcal{S}_G$ of separable $\overline{\ell}$-algebras we define a quiver
 $Q_0(\mathcal{S}_G)$ as follows 
 \begin{itemize}
 \item{Vertices :  for any vertex $v \in V$ of $G$ take $n_v$ vertices $\{ \mu^{(v)}_1,\hdots,\mu^{v}_{n_v} \}$. }
 \item{Arrows : fix an orientation $\vec{G}$ on all of the edges of $G$. For any edge $\xymatrix{\vtx{v} \ar@{-}[r]^e & \vtx{w} }$ in $G$ we add for each $1 \leq i \leq n_v$ and each $1 \leq j \leq n_w$ precisely $n^{(e)}_{ij}$ arrows between the vertices $\mu^{(v)}_i$ and $\mu^{(w)}_j$ oriented in the same way as the edge $e$ in $\vec{G}$.}
 \end{itemize}
 We call $Q_0(\mathcal{S}_G)$ the {\em Zariski quiver} of the graph of separable algebras $\mathcal{S}_G$.
 \end{definition}
 
 The {\em representation space} $rep_{\alpha}~Q_0(\mathcal{S}_G)$ is the affine $\overline{\ell}$-space
 \[
 rep_{\alpha}~Q_0(\mathcal{S}_G) = \bigoplus_{\xymatrix{\vtx{v} \ar[r]^{e} & \vtx{w}}} \oplus_{i=1}^{n_v} \oplus_{j=1}^{n_w} M_{\alpha^{(w)}_j \times \alpha^{(v)}_i}(\overline{\ell}) \]
 and two $\alpha$-dimensional representations are said to be {\em isomorphic} if they are conjugated via the natural base-change action of $GL(\alpha) = \times_{v \in V} \times_{i=1}^n GL(\alpha^{(v)}_i)$.

A dimension vector $\alpha = (\alpha^{(v)}_i~:~v \in V, 1 \leq i \leq n_v)$ for $Q_0(\mathcal{S}_G)$ is said to be an {\em $n$-dimension vector} if  the following numerical conditions are satisfied
 \[
 \sum_{i=1}^{n_v} d_v(i) \alpha^{(v)}_i = n  \]
for all $v \in V$.

 For any edge $\xymatrix{\vtx{v} \ar[r]^{e} & \vtx{w}}$ we denote by $Q_e$ the {\em bipartite} subquiver of $Q_0(\mathcal{S}_G)$ on the vertices $\{ \mu_1^{(v)},\hdots,\mu_{n_v}^{(v)} \}$, $\{ \mu_1^{(w)},\hdots,\mu_{n_w}^{(w)} \}$ and the $n_{ij}^{(e)}$ arrows between $\mu_i^{(v)}$ and $\mu_j^{(w)}$ determined by the embeddings $S_e \rInto S_v$ and $S_e \rInto S_w$.
 
 \begin{definition}
 Let $\alpha$ be an $n$-dimension vector, $M \in rep_{\alpha}~Q_0(\mathcal{S}_G)$ and $e \in E$ :
 \begin{itemize}
 \item{$M$ is said to be $e$-semistable iff for all $Q_e$- subrepresentations $N$ of $M | Q_e$ of dimension vector $(n_1,\hdots,n_{n_v},n'_1,\hdots,n'_{n_w})$ we have
 \[
 \sum_{i=1}^{n_w} n'_i d_w(i) \geq \sum_{i=1}^{n_v} n_i d_v(i) \]}
 \item{$M$ is said to be $e$-stable iff for all proper $Q_e$-subrepresentations $N$  of $M | Q_e$ of
 dimension vector $(n_1,\hdots,n_{n_v},n'_1,\hdots,n'_{n_w})$ we have
 \[
 \sum_{i=1}^{n_w} n'_i d_w(i) > \sum_{i=1}^{n_v} n_i d_v(i) \]}
 \item{$M$ is said to be $\mathcal{S}_G$-semistable (resp. $\mathcal{S}_G$-stable) iff $M$ is $e$-semistable (resp. $e$-stable) for all edges $e \in E$.}
 \end{itemize}
 \end{definition}
 
 The relevance of the quiver $Q_0(\mathcal{S}_G)$ and the introduced terminology is contained in the next result.
 
  \begin{proposition} Every $n$-dimensional representation $\pi_1(\mathcal{S}_G) \rTo^{\phi} M_n(\overline{\ell})$ determines (and is determined by) an $\mathcal{S}_G$-semistable representation $M_{\phi} \in rep_{\alpha}~Q_0(\mathcal{S}_G)$ for some $n$-dimension vector $\alpha$. Moreover, if $\phi$ and $\phi'$ are isomorphic representations of $\pi_1(\mathcal{S}_G)$, then $M_{\phi}$ and $M_{\phi'}$ are isomorphic as quiver representations.
 \end{proposition}
 
 \begin{proof} Let $N = \overline{\ell}^{n}_{\phi}$ be the $n$-dimensional module determined by $\phi$. For each vertex $v \in V$ we have a decomposition by restricting $N$ to the separable subalgebra $S_v$
 \[
 N \downarrow_{S_v} \simeq V_{1,v}^{\oplus \alpha_1^{(v)}} \oplus \hdots \oplus V_{n_v,v}^{\oplus \alpha_{n_v}^{(v)}} \]
 where the $V_{i,v}$ are the distinct simple modules of $S_v$ of dimension $d_v(i)$. Choose an $\overline{\ell}$-basis $\mathcal{B}_v$ of $N \downarrow_{S_v}$ compatible with this decomposition. These decompositions determine an $n$-dimension vector $\alpha$.  For any edge $\xymatrix{\vtx{v} \ar[r]^{e} & \vtx{w}}$ the embeddings $S_e \rInto^{\alpha} S_v$ and $S_e \rInto^{\beta} S_w$ determine two $n$-dimensional $S_e$-representations
 \[
 (N \downarrow_{S_v}) \downarrow_{S_e}^{\alpha} \qquad \text{and} \qquad (N \downarrow_{S_w}) \downarrow_{S_e}^{\beta} \]
 which, by construction of $\pi_1(\mathcal{S}_G)$ are isomorphic. That is, the basechange map $\mathcal{B}_v \rTo^{\psi_{vw}} \mathcal{B}_w$ is an invertible element of
 \[
 Hom_{S_e}(N \downarrow_{S_v}, N \downarrow_{S_w}) = \oplus_{i=1}^{n_v} \oplus_{j=1}^{n_w} M_{\alpha_j^{(w)} \times \alpha_i^{(v)}}(Hom_{S_e}(V_{i,v},V_{j,w})) \]
and hence $\psi_{vw}$ determines a representation of the bipartite quiver $Q_e$ of dimension vector $\alpha | Q_e$. Repeating this for all edges $e \in E$ we obtain a representation $M_{\phi} \in rep_{\alpha}~Q_0(\mathcal{S}_G)$. Invertibility of the map $\psi_{vw}$ is equivalent to $M_{\phi}$ being $e$-semistable, so $M_{\phi}$ is $\mathcal{S}_G$-semistable. Isomorphic representations $\phi$ and $\phi'$ determine isomorphic vertex-decompositions whence, by Schur's lemma, bases which are transferred into each other via an element of $GL(\alpha)$ and hence the quiver representations $M_{\phi}$ and $M_{\phi'}$ are isomorphic. From the construction of the fundamental algebra $\pi_1(\mathcal{S}_G)$ it follows that one can reverse this procedure to construct on $n$-dimensional representation of $\pi_1(\mathcal{S}_G)$ from a $\mathcal{S}_G$-stable representation $M \in rep_{\alpha}~Q_0(\mathcal{S}_G)$ for some $n$-dimension vector $\alpha$.
 \end{proof}
 
 Under this correspondence simple $\pi_1(\mathcal{S}_G)$-representations correspond to $\mathcal{S}_G$-stable representations. If $\alpha$ is an $n$-dimension vector such that $rep_{\alpha}~Q_0(\mathcal{S}_G)$ contains $\mathcal{S}_G$-stable representations (which then form a Zariski open subset), then $\alpha$ is a {\em Schur root} of $Q_0(\mathcal{S}_G)$ and consequently the dimension of the classifying variety is equal to
$ 1 - \chi_0(\alpha,\alpha)$ where $\chi$ is the {\em Euler form} of the quiver $Q_0(\mathcal{S}_G)$. For this result and related material on Schur roots we refer to \cite{Schofieldgeneric}. 

\begin{proposition} Isomorphism classes of simple $n$-dimensional representations of $\pi_1(\mathcal{S}_G)$ are parametrized by the points of a smooth quasi-affine variety (possibly with several irreducible components)
\[
isosimp_n~\pi_1(\mathcal{S}_G) = \bigsqcup_{\alpha} isosimp_{\alpha}~\pi_1(\mathcal{S}_G) \]
where $\alpha$ runs over all $n$-dimension vectors such that $rep_{\alpha}~Q_0(\mathcal{S}_G)$ contains $\mathcal{S}_G$-stable representations. These components have dimensions
\[
dim~isosimp_{\alpha}~\pi_1(\mathcal{S}_G) = 1 - \chi_0(\alpha,\alpha) \]
where $\chi_0$ is the Euler form of the quiver $Q_0(\mathcal{S}_G)$.
\end{proposition}

As an example consider the modular group $SL_2(\Z)$ which is the amalgamated product $\Z_4 \ast_{\Z_2} \Z_6$, see for example \cite[I \S 7]{Dicks}. If $char(\overline{\ell}) \not= 2,3$ the group-algebra $\overline{\ell} SL_2(\Z)$ is the fundamental algebra of the graph of separable $\overline{\ell}$-algebras
\[
\xymatrix{\vtx{v} \ar[r]^e & \vtx{w}} \qquad \text{with} \qquad S_v = \overline{\ell} \Z_4~\quad~S_w = \overline{\ell} \Z_6~\quad~S_e = \overline{\ell} \Z_2 \]
As all simples are one-dimensional (determined by their eigenvalue), it is easy to verify that the zero quiver $Q_0(\overline{\ell} SL_2(\Z))$ has the following form
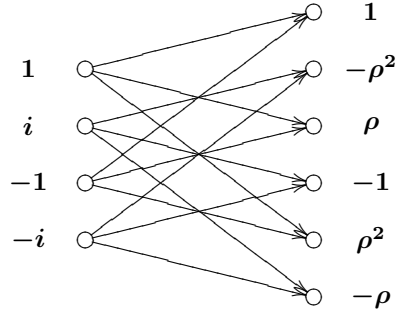
\[
\xy /r.18pc/:
\POS (40,0) *\cir<3pt>{}="b1",
 (40,10) *\cir<3pt>{}="b2",
 (40,20) *\cir<3pt>{}="b3",
 (40,30) *\cir<3pt>{}="b4",
 (40,40) *\cir<3pt>{}="b5",
 (40,50) *\cir<3pt>{}="b6",
 (0,10) *\cir<3pt>{}="a1",
 (0,20) *\cir<3pt>{}="a2",
 (0,30) *\cir<3pt>{}="a3",
 (0,40) *\cir<3pt>{}="a4",
 (-10,10) *\txt{$-i$},
 (-10,20) *\txt{$-1$},
 (-10,30) *\txt{$i$},
 (-10,40) *\txt{$1$},
 (50,0) *\txt{$-\rho$},
 (50,10) *\txt{$\rho^2$},
 (50,20) *\txt{$-1$},
 (50,30) *\txt{$\rho$},
 (50,40) *\txt{$-\rho^2$},
 (50,50) *\txt{$1$}
 \POS "a1" \ar "b1"
  \POS "a1" \ar "b3"
    \POS "a1" \ar "b5"
  \POS "a2" \ar "b2"
  \POS "a2" \ar "b4"
  \POS "a2" \ar "b6"
  \POS "a3" \ar "b1"
  \POS "a3" \ar "b3"
  \POS "a3" \ar "b5"
  \POS "a4" \ar "b2"
  \POS "a4" \ar "b4"
  \POS "a4" \ar "b6"
 \endxy
 \]
($\rho$ is a primitive $3$rd root of unity) which is the disjoint union of two copies of the quiver associated to $PSL_2(\Z)$ in \cite{Westbury}.

The congruence subgroup $\Gamma_0(2) = \{ \begin{bmatrix} a & b \\ c & d \end{bmatrix} \in SL_2(\Z)~\text{ with $c$ even}~\}$ is the fundamental group of the graph of finite groups
\[
\xymatrix{\vtx{v} \ar@{-}[r]^e & \vtx{w} \ar@{-}@(ur,dr)^f}~\qquad G_w = G_e = G_f = \Z_2,~G_v = \Z_4 
\]
If $char(\overline{\ell}) \not= 2$, the group algebra $\overline{\ell} \Gamma_0(2)$ is the fundamental algebra of a graph of separable $\overline{\ell}$-algebras and the zero quiver $Q_0(\overline{\ell} \Gamma_0(2))$ has the following form
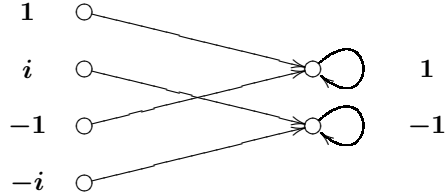
\[
\xy /r.18pc/:
\POS 
 (40,20) *\cir<3pt>{}="b3",
 (40,30) *\cir<3pt>{}="b4",
  (0,10) *\cir<3pt>{}="a1",
 (0,20) *\cir<3pt>{}="a2",
 (0,30) *\cir<3pt>{}="a3",
 (0,40) *\cir<3pt>{}="a4",
 (-10,10) *\txt{$-i$},
 (-10,20) *\txt{$-1$},
 (-10,30) *\txt{$i$},
 (-10,40) *\txt{$1$},
  (60,20) *\txt{$-1$},
 (60,30) *\txt{$1$},
  \POS "a1" \ar "b3"
   \POS "a2" \ar "b4"
  \POS "a3" \ar "b3"
   \POS "a4" \ar "b4"
   \POS "b3" \ar@(ur,dr) "b3"
   \POS "b4" \ar@(ur,dr) "b4"
 \endxy
 \]

\begin{definition} For a graph $\mathcal{S}_G$ of separable $\overline{\ell}$-algebras we define a quiver $Q_1(\mathcal{S}_G)$ as follows
\begin{itemize}
\item{Vertices : Let $\{ \alpha_1,\hdots,\alpha_k \}$ be the minimal set of generators for the sub-semigroup of dimension vectors $\alpha$ of $Q_0(\mathcal{S}_G)$ which are $n$-dimension vectors for some $n \in \N$ and such that $rep_{\alpha}~Q_0(\mathcal{S}_G)$ contains $\mathcal{S}_G$-semistable representations. The vertices $\{ \nu_1,\hdots,\nu_k \}$ are in one-to-one correspondence with these generators $\{ \alpha_1,\hdots,\alpha_k \}$.}
\item{Arrows : The number of directed arrows in $Q_1(\mathcal{S}_G)$ from $\nu_i$ to $\nu_j$
\[
\#~\{~\xymatrix{\vtx{i} \ar[r] & \vtx{j}}~\} = \delta_{ij} - \chi_0(\alpha_i,\alpha_j) \]
where $\chi_0$ is the Euler-form of the Zariski quiver $Q_0(\mathcal{S}_G)$.}
\end{itemize}
We call $Q_1(\mathcal{S}_G)$ the one-quiver of the graph of separable algebras $\mathcal{S}_G$.
\end{definition}

The one-quiver $Q_1(\mathcal{S}_G)$ allows us to determine the components $rep_{\alpha}~\pi_1(\mathcal{S}_G)$ which contain (a Zariski open subset of) simple representations. Remark that the description of Schur roots is a lot harder than that of dimension vectors of simple representations.

\begin{proposition} If $\alpha = c_1 \alpha_1 + \hdots + c_k \alpha_k \in comp~\pi_1(\mathcal{S}_G)$ then the component $rep_{\alpha}~\pi_1(\mathcal{S}_G)$ contains simple representations if and only if
\begin{itemize}
\item{
\[
\chi_1(\gamma,\epsilon_i) \leq 0 \qquad \text{and} \qquad \chi_1(\epsilon_i,\gamma) \leq 0 \]
for all $1 \leq i \leq k$ where $\gamma = (c_1,\hdots,c_k)$ and $\epsilon_i = (\delta_{1i},\hdots,\delta_{ki})$ and where $\chi_1$ is the Euler form of the one quiver $Q_1(\mathcal{S}_G)$}
\item{$supp(\gamma)$ is a strongly connected subquiver of $\pi_1(\mathcal{S}_G)$ and if $supp(\gamma)$ is of extended Dynkin type $\tilde{A}_l$ then all non-zero components of $\gamma$ must be equal to one.}
\end{itemize}
\end{proposition}

\begin{proof} Follows from the proof of proposition~\ref{simpele}.
\end{proof}

If $char(\overline{\ell}) = 0$ one can apply Luna slice machinery to construct a Zariski open subset of all simple representations in $rep_{\alpha}~\pi_1(\mathcal{S}_G)$ from the knowledge of low-dimensional simples. For example, suppose we have found simple representations
\[
S_i \in rep_{\alpha_i}~\pi_1(\mathcal{S}_G) \qquad \text{for all $1 \leq i \leq k$} \]
and consider the point $M$ in the affine space $rep_{\alpha}~Q_0(\mathcal{S}_G)$ determined by the semi-simple representation of $\pi_1(\mathcal{S}_G)$
\[
M = S_1^{\oplus c_1} \oplus \hdots \oplus S_k^{\oplus c_k} \]
then the normal space to the $GL(\alpha)$-orbit $\Oscr(M)$ is isomorphic to $Ext^1_{\pi_1(\mathcal{S}_G)}(M,M)$ which we have seen can be identified to $rep_{\gamma}~Q_1(\mathcal{S}_G)$. 

\begin{proposition} Let $\alpha = c_1 \alpha_1 + \hdots + c_k \alpha_k$ be a component such that $rep_{\alpha}~\pi_1(\mathcal{S}_G)$ contains simple representations. In the affine space $rep_{\alpha} Q_0(\mathcal{S}_G)$ identify the normal space to the orbit $\Oscr(M)$ of the semi-simple representation $M$ (as above) with
\[
N_M = \{ M + V~|~V \in rep_{\gamma}~Q_1(\mathcal{S}_G)~\} \]
where $\gamma = (c_1,\hdots,c_k)$.
Then, $GL(\alpha).N_M$ contains a Zariski open subset of all $\alpha$-dimensional simple representations of $\pi_1(\mathcal{S}_G)$.
\end{proposition}

\begin{proof} This is a special case of the Luna slice result applied to the local quiver setting.
In fact, one can generalize this result to other known semi-simple representations $N$ of $\pi_1(\mathcal{S}_G)$ but then one has to replace $Q_1(\mathcal{S}_G)$ by the {\em local quiver} $Q_N$ of $N$.
\end{proof}

In the $SL_2(\Z)$ example, $comp(\overline{\ell} SL_2(\Z))$ is generated by the $12$ components of two-dimensional representations of $Q_0(\overline{\ell} SL_2(\Z))$
\[
\nu_{ij} = (\delta_{1i},\hdots,\delta_{4i};\delta_{1j},\hdots,\delta_{6j}) \qquad 1 \leq i \leq 4, 1 \leq j \leq 6
\]
for which $i$ and $j$ are both even or both odd.
From this the structure of the one quiver $Q_1(\overline{\ell} SL_2(\Z))$ (corresponding to the $12$ one-dimensional simples of $\overline{\ell} SL_2(\Z)$) can be verified to be
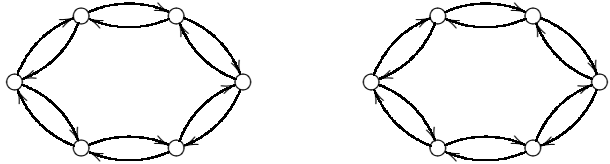
\[
\xy /.15pc/:
\POS (0,0) *\cir<3pt>{} ="S11",
\POS (14,14) *\cir<3pt>{} ="S23",
\POS (-14,14) *\cir<3pt>{} ="S22",
\POS (14,34) *\cir<3pt>{} ="S12",
\POS (-14,34) *\cir<3pt>{} ="S13",
\POS (0,48) *\cir<3pt>{}="S21",
\POS"S11" \ar@/^1ex/ "S23"
\POS"S23" \ar@/^1ex/ "S11"
\POS"S12" \ar@/^1ex/ "S23"
\POS"S23" \ar@/^1ex/ "S12"
\POS"S12" \ar@/^1ex/ "S21"
\POS"S21" \ar@/^1ex/ "S12"
\POS"S21" \ar@/^1ex/ "S13"
\POS"S13" \ar@/^1ex/ "S21"
\POS"S13" \ar@/^1ex/ "S22"
\POS"S22" \ar@/^1ex/ "S13"
\POS"S11" \ar@/^1ex/ "S22"
\POS"S22" \ar@/^1ex/ "S11"
\endxy
\qquad~\qquad
\xy /.15pc/:
\POS (0,0) *\cir<3pt>{} ="S11",
\POS (14,14) *\cir<3pt>{} ="S23",
\POS (-14,14) *\cir<3pt>{} ="S22",
\POS (14,34) *\cir<3pt>{} ="S12",
\POS (-14,34) *\cir<3pt>{} ="S13",
\POS (0,48) *\cir<3pt>{}="S21",
\POS"S11" \ar@/^1ex/ "S23"
\POS"S23" \ar@/^1ex/ "S11"
\POS"S12" \ar@/^1ex/ "S23"
\POS"S23" \ar@/^1ex/ "S12"
\POS"S12" \ar@/^1ex/ "S21"
\POS"S21" \ar@/^1ex/ "S12"
\POS"S21" \ar@/^1ex/ "S13"
\POS"S13" \ar@/^1ex/ "S21"
\POS"S13" \ar@/^1ex/ "S22"
\POS"S22" \ar@/^1ex/ "S13"
\POS"S11" \ar@/^1ex/ "S22"
\POS"S22" \ar@/^1ex/ "S11"
\endxy
\]
Here, the vertices of the first component correspond (in cyclic order) to $\nu_{11},\nu_{33},\nu_{15},\nu_{31},\nu_{13},\nu_{35}$ and those of the second component (in cyclic order) to $\nu_{22},\nu_{44},\nu_{26},\nu_{42},\nu_{24},\nu_{46}$. Applications to the representation theory of the modular group $SL_2(\Z)$ and its central extension $B_3$ (the third braid group) will be given elsewhere.

In the $\Gamma_0(2)$ example, $comp(\overline{\ell} \Gamma_0(2))$ is generated by the $4$ dimension vectors
\[
(1,0,0,0;1,0),~(0,0,1,0;1,0),~(0,1,0,0;0,1),~(0,0,0,1;0,1) \]
and one verifies that the one-quiver $Q_1(\overline{\ell} \Gamma_0(2))$ has the following form
\vskip 3mm
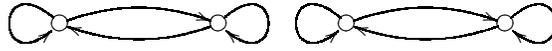
\[
\xymatrix{ \vtx{} \ar@/^/[rr] \ar@(ul,dl) & & \vtx{} \ar@/^/[ll] \ar@(ur,dr)} \qquad~\qquad
\xymatrix{ \vtx{} \ar@/^/[rr] \ar@(ul,dl) & & \vtx{} \ar@/^/[ll] \ar@(ur,dr)} 
\]

\section*{Appendix : The component coalgebra $coco(A)$}

Over an algebraically closed field $\overline{\ell}$ we have seen that the component semigroup and Euler form contain useful information on the finite dimensional representations of an $\overline{\ell}$-qurve. Clearly, one can repeat all arguments verbatim for an arbitrary $\ell$ by restricting at those components which contain $\ell$-rational points. However, this sub-semigroup $\wis{comp}(A)$ of $\wis{comp}(A \otimes \overline{\ell})$ is usually too small to be of interest.

\begin{example} Let $\ell \subset L$ be a finite separable field extension of dimension $k$. As $L$ is a simple algebra, all its finite dimensional representations are of the form $L^{\oplus a}$ and hence only components of $rep_n~L$ containing $\ell$-rational points exist when $k | n$. Over the algebraic closure we have
\[
L \otimes \overline{\ell} = \underbrace{\overline{\ell} \times \hdots \times \overline{\ell}}_k \]
whence $comp(L \otimes \overline{\ell}) \simeq \N^k$ generated by the factors of $L \otimes \overline{\ell}$. We have $comp(A) \subset comp(A \otimes \overline{\ell})$ sending the generator $k$ to $(1,\hdots,1)$.
\end{example}

We recall some standard facts from \cite[Chp. 1]{Iversen} on unramified commutative algebras over an arbitrary basefield $\ell$.
A commutative affine $\ell$-algebra $C$ is said to be {\em unramified} whenever
\[
C \otimes \overline{\ell} \simeq \overline{\ell} \times \hdots \times \overline{\ell} \]
It is well known that all unramified $\ell$-algebras are of the form
\[
C \simeq L_1 \times \hdots \times L_k \]
where each $L_i$ is a finite dimensional separable field extension of $\ell$. From this it follows that subalgebras, tensorproducts and epimorphic images of unramified $\ell$-algebras are again unramified. As a consequence, an affine commutative $\ell$-algebra $C$ has a unique {\em maximal unramified $\ell$-subalgebra} $\pi_0(C)$. In case $C = \ell[X]$ is the coordinate algebra of an affine $\ell$-scheme $X$, the algebra $\pi_0(C)$ contains all information about the connected components of $X$.
Recall that an affine $\ell$-scheme $X$ (or its coordinate algebra $\ell[X]$) is said to be {\em connected} if $\ell[X]$ contains no non-trivial idempotents and is called {\em geometrically connected} if $\ell[X] \otimes \overline{\ell}$ is connected. We summarize \cite[I.7]{Iversen} in

\begin{proposition} For an affine $\ell$-scheme $X$ we have
\begin{enumerate}
\item{$X$ is connected iff $\pi_0(\ell[X])$ is a field.}
\item{$X$ is geometrically connected iff $\pi_0(\ell[X]) = \ell$.}
\item{If $X$ is connected and has an $\ell$-rational point, then $X$ is geometrically connected.}
\item{If $\pi_0(\ell[X]) = L_1 \times \hdots \times L_k$ with all $L_i$ separable field extensions of $\ell$, then $X$ has exactly $k$ connected components.}
\item{If $Y$ is an affine $\ell$-scheme and $X \rTo Y$ a morphism, then $\pi_0(\ell[Y]) \rTo \pi_0(\ell[X])$ is an $\ell$-algebra morphism.}
\item{If $Y$ is an affine $\ell$-scheme, then the natural map
\[
\pi_0(\ell[X]) \otimes \pi_0(\ell[Y]) \rTo \pi_0(\ell[X] \otimes \ell[Y]) = \pi_0(\ell[X \times Y]) \]
is an $\ell$-algebra isomorphism.}
\end{enumerate}
\end{proposition}

\begin{definition}
For $A$ an $\ell$-qurve consider the sum-maps
\[
rep_n~A \times rep_m~A \rTo rep_{m+n}~A \]
which determine $\ell$-algebra morphisms
\[
\Delta_{m,n}~:~\pi_0(\ell[rep_{m+n}~A]) \rTo \pi_0(\ell[rep_n~A]) \otimes \pi_0(\ell[rep_m~A]) \]
Denote $\pi_0(n) = \pi_0(\ell[rep_n~A])$ and consider the graded $\ell$-vectorspace
\[
coco(A) = \pi_0(0) \oplus \pi_0(1) \oplus \pi_0(2) \oplus \hdots 
\]
Define a coalgebra structure by taking as the {\em comultiplication} map 
\[
coco(A) \rTo^{\Delta} coco(A) \otimes coco(A)
\]
\[
\sum_{m+n=N} \Delta_{m,n}~:~ \pi_0(N) \rTo \sum_{n+m=N} \pi_0(n) \otimes \pi_0(m) \]
and as the {\em counit} $\wis{coco}(A) \rOnto^{\epsilon} \pi_0(0) = \ell$.
We call $(coco(A),\Delta,\epsilon)$ the {\em component coalgebra} of the $\ell$-qurve $A$.
\end{definition}

In fact, it follows from the foregoing proposition that $coco(A)$ is in fact a {\em mock bialgebra}, that is a bialgebra without a unit-map. Recall that if $G$ is a finite group, its {\em function bialgebra} $func(G)$ is the space of all $\ell$-valued functions on $G$ with pointwise multiplication and co-multiplication induced by
\[
\Delta(x_g) = \sum_{g'.g"=g} x_{g'} \otimes x_{g"} \]
where $x_h$ is the function mapping $h \mapsto 1$ and all other $h' \in G$ to zero. If $G$ is no longer finite, $func(G)$ is still a mock bialgebra.

\begin{proposition} If $A$ is an $\ell$-qurve, then there is an isomorphism of mock bialgebras
\[
coco(A) \otimes \overline{\ell} \simeq func(comp(A \otimes \overline{\ell})) \]
and hence $coco(A)$ contains enough information to reconstruct the component semigroup $comp(A \otimes \overline{\ell})$. Alternatively, the Galois group $Gal(\overline{\ell}/\ell)$ acts on $A \otimes \overline{\ell}$ and hence on $comp(A \otimes \overline{\ell})$ and the function coalgebra. The component coalgebra $coco(A)$ can be obtained by Galois descent
\[
coco(A) = func(comp(A \otimes \overline{\ell}))^{Gal(\overline{\ell}/\ell)} \]
\end{proposition}

 \end{document}